\documentclass[a4paper,12pt]{article}
                                                                                
\usepackage{amssymb}

\def\mmp{{{\cal M}_{\lambda,\Lambda}^+}}
\def\mmm{{{\cal M}_{\lambda,\Lambda}^-}}
\def\mmpm{{{\cal M}_{\lambda,\Lambda}^\pm}}

\newcommand{\RR}{\hbox{\bf R}}

\newcommand{\be}{\begin{equation}}
\newcommand{\ee}{\end{equation}}
\newcommand{\bea}{\begin{eqnarray}}
\newcommand{\eea}{\end{eqnarray}}

\newcommand{\equ}[1]{(\ref{#1})}

\newcommand{\R}{{ I\!\!R}}
\newcommand{\N}{{ I\!\!N}}
 
\def\sobre#1#2{\lower 1ex \hbox{ $#1 \atop #2 $ } }
 
\def\bajo#1#2{\raise 1ex \hbox{ $#1 \atop #2 $ } }

\newtheorem{theorem}{Theorem}[section]
\newtheorem{proposition}{Proposition}[section]
\newtheorem{lemma}{Lemma}[section]

\newtheorem{remark}{Remark}[section]
\newcommand{\bremark}{\begin{remark} \em}
\newcommand{\eremark}{\end{remark} }

\begin{document}

\centerline{\bf LARGE CRITICAL EXPONENTS FOR SOME} 
\centerline{\bf  SECOND ORDER UNIFORMLY 
ELLIPTIC OPERATORS}
\bigskip
\bigskip
\small{
\centerline{ by}
\centerline{Maria J. ESTEBAN${}^a$,
Patricio L. FELMER${}^{b}$}
\centerline{and}
\centerline{Alexander QUAAS${}^{c}$}
  }                                                                              
\bigskip
\small{
\centerline{${}^{a}$ Ceremade UMR CNRS 7534, Universit\'e Paris-Dauphine}
\centerline{75775 Paris Cedex 16, FRANCE.} 

\medskip

\centerline{${}^{b}$ Departamento de Ingenier\'{\i}a  Matem\'atica, }
\centerline{ and Centro de Modelamiento Matem\'atico, UMR2071 CNRS-UChile}
\centerline{ Universidad de Chile, Casilla 170 Correo 3,}
 \centerline{ Santiago, CHILE.}
       
\medskip

                                                                         \centerline{${}^{c}$ Departamento de  Matem\'atica, Universidad Santa Mar\'{\i}a,}
\centerline{Casilla: V-110, Avda. Espa\~na 1680, Valpara\'{\i}so, CHILE.} 
}

\bigskip
\bigskip
\bigskip
\bigskip
\bigskip

\section{Introduction}

Associated to the Laplacian we have the Sobolev critical exponent, which is the largest
 number  $p^*_N$ having the property  that the semi-linear equation 
\begin{eqnarray}\label{generaleqI}
\Delta u +u^p&=&0,\quad \mbox{in}\quad \Omega\label{1I}\\
u&=&0,\quad \mbox{on}\quad \partial\Omega\label{2I},
\end{eqnarray}
possesses a positive solution whenever $1<p<p^*_N$ and for any given bounded domain
$\Omega$. This number, which depends on the dimension $N$, is given by 
$p^*_N=(N+2)/(N-2)$ and its name derives from the fact that the Sobolev space 
$H^1_0(\Omega)$ embeds itself continuously into $L^q(\Omega)$,
 for all domain $\Omega$ of $\R^N$, if and only if $1<q\leq p^*_N+1$, and 
the embedding is compact if and only if $1<q< p^*_N+1$ (and if the domain is bounded or at least very small at infinity).
Moreover, if $p\ge p^*_N$ and the domain is star-shaped, Pohozaev's identity implies that the above equation does not have a positive solution 
(see \cite{P}). The Sobolev exponent has a
 dual property if the domain is $\R^N$, actually the equation has a positive solution whenever $p\ge p^*_N$ and it does not have a solution if 
$1<p<p^*_N$. 

If we replace the Laplacian by 
any linear second order uniformly elliptic operator with $C^1$ coefficients, 
say $Lu=\sum_i\sum_j a_{ij} \frac{\partial u^2}{\partial x_i\partial x_{j}}$ with $a_{ij}\in C^1$, then 
the semi-linear problem 
\begin{eqnarray}
L u +u^p&=&0,\quad \mbox{in}\quad \Omega\label{1}\\
u&=&0,\quad \mbox{on}\quad \partial\Omega\label{2},
\end{eqnarray}
has a positive solution for the same range of values of $p$, namely $1<p<p^*_N$.
That is, the existence property of the Sobolev exponent remains valid 
for all  operators in this class.

In this note we  consider two classes of uniformly second order 
elliptic operators for which the critical exponents are drastically changed
with respect to $p_N^*$, in the case of radially symmetric solutions. Our
aim
is to prove 
that the corresponding existence property for these critical exponents persists when the
domain is perturbed.

Our first class corresponds to the so-called Pucci's extremal operators 
\cite{CPucci}, \cite{CPucci1} and \cite{cc}. Given positive numbers 
$0<\lambda\le \Lambda$ we consider the operator
$
{\cal M}_{\lambda,\Lambda}^+(D^2u)\label{P}
$
where for any $N\times N$ symmetric matrix $M$,  $${\cal M}_{\lambda,\Lambda}^+(M)=\Lambda \sum_{e_i>0}e_i+\lambda  
\sum_{e_i<0}e_i,
$$
  $e_i=e_i(M)$ being $M$'s eigenvalues. The case of the operator 
${\cal M}_{\lambda,\Lambda}^-(D^2u)$ is also considered and it is defined
by exchanging the roles between $\lambda$ and $\Lambda$ above. 

Pucci's extremal operators appear in the context of stochastic
control when the diffusion coefficient is a control variable, see
the book of Bensoussan and J.L. Lions \cite{ben} or the papers of
P.L. Lions \cite{lionspI}, \cite{lionspII}, \cite{lionspIII} for the relation between a general Hamilton-Jacobi-Bellman
and stochastic control.
They also provide natural extremal equations in the sense that if $F$
 is any (from linear to fully nonlinear) uniformly elliptic operator, with ellipticity
constants $\lambda$, $\Lambda$, and depends only on the Hessian
$D^2u$, then 

\be\label{ineg_extremale} \mmm(M) \leq F(M) \leq \mmp(M) \ee for
any symmetric matrix $M$.  Moreover, these operators are also extremal with respect to
the first half eigenvalue of all second order elliptic operators with constant coefficients and ellipticity constants between  $\lambda$ and $\Lambda$ (see for instance \cite{busca-esteban-quaas}).

It is obvious that when $\lambda=\Lambda$, then $\mmpm$ coincides with 
a multiple of the Laplace
operator. We also notice that given any number $s\in  [\lambda,\Lambda]$ the operator $s\Delta$ belongs to the class defined by \equ{ineg_extremale}.

The second family of operators that we consider are 
 defined as
\be
Q^+_{\lambda,\Lambda}u=\lambda\Delta u +(\Lambda-\lambda) Q^0u,\label{R}
\ee
where  $Q^0$ is the second order operator
$$
Q^0u=\sum_{i=1}^N\sum_{j=1}^N\frac{x_ix_j}{|x|^2}\frac{\partial^2 u}
{\partial x_i\partial x_j}
\,.$$
These operators are also considered by Pucci \cite{CPucci}, being 
extremal with respect to some spectral properties. 
We notice that these operators  belong to the class defined by 
\equ{ineg_extremale} 
and when $\lambda=\Lambda$ they also become a multiple of the Laplacian. 
If we interchange the role of $\lambda$ by $\Lambda$
in definition \equ{R}, then we obtain the operator $Q^-_{\lambda,\Lambda}$, which is also considered later.

The operators $\,{\cal M}_{\lambda,\Lambda}^\pm \,$ are autonomous, but not linear, 
even if they enjoy some properties of the Laplacian. The operators $\,Q_{\lambda,\Lambda}^\pm \,$ 
are still linear, but their coefficients are not continuous at 
the origin.
In both cases, when one considers a ball  and the set of radially symmetric functions in it, there are  critical exponents for the operators $\,{\cal M}^+ \,$ and $\,Q^+\,$ 
which are greater than the Sobolev exponent $p^*_N$. On the contrary, for the 
operators $Q_{\lambda,\Lambda}^-$ and $\mmm$, the critical exponents for radially symmetric solutions in a ball are smaller than the Sobolev exponent $p^*_N$.
These facts where proved  in \cite{felmerquaas, felmerquaas1} for $\,{\cal M}_{\lambda,\Lambda}^\pm \,$
 and for $Q_{\lambda,\Lambda}^+$ and $Q_{\lambda,\Lambda}^-$ the proof is given here, in Section 2.

More precisely, in the case of operators $Q^+_{\lambda,\Lambda}$, 
there exists a number 
$\widetilde N_+=\frac{\lambda}{\Lambda}(N-1)+1\,,$
such that if $\Omega$ is a ball of $\R^N$ and if $\widetilde N_+>2$,  
\equ{1}-\equ{2} 
has a unique positive radially symmetric solution for 
any $1<p<(\widetilde N_++2)/(\widetilde N_+-2)$ and no positive 
radially symmetric solution for 
$p\geq (\widetilde N_++2)/(\widetilde N_+-2)$. Notice that for any $\lambda<\Lambda$,
$\widetilde N_+<N$ and so, the critical exponent here is strictly larger than the Sobolev critical exponent $p^*_N$.

In the case of the Pucci's extremal operators $\mmp$, the critical exponent is   a number $p^*_+$ such that 
 $$p^*_N<p^*_+< \frac{\widetilde N_++2}{\widetilde N_+-2}\,.$$
The number $p^*_+$ depends on $\lambda,$ $\Lambda$ and the dimension $N$, 
however an explicit formula for it is
 not known. 

Similarly, for the operator $Q^-_{\lambda,\Lambda}$ we may also define a dimension like number
$\widetilde N_-=\frac{\Lambda}{\lambda}(N-1)+1\,$, so that
 its critical exponent in the radially symmetric case is precisely
$  ({\widetilde N_-+2})/({\widetilde N_--2}).$
In the case of the operator $\mmm$ we recall that the critical exponent of the operator $\mmm$ is a number $p^*_-$
satisfying
$$
\frac{\widetilde N_-+2}{\widetilde N_--2}\,<p^*_-< p_N^*,$$
as it was shown in \cite{felmerquaas}.

It is the purpose of this note to prove that this phenomenom of
critical  exponent increase (or decrease) does not appear
only in the radially symmetric case.
 By a perturbation argument, based on a work by Dancer \cite{dancer}, we 
show that these critical exponents, with respect to 
existence properties in bounded domains,
persist when the ball is perturbed not necessarily in a radial manner.  This result provides us with evidence that the critical exponents  for these
operators, obtained in radial versions,  are also the critical exponents in the general case.

At this point we would like to stress some surprising properties of  the
critical exponents of operators in the class given by \equ{ineg_extremale}.
For the first property we consider all linear elliptic operators with bounded coefficients
and belonging to the class defined by \equ{ineg_extremale}. If we take
the $L^\infty$ topology for the coefficients of these operators, we 
see that the critical exponent
is not a continuous function of the operator. In particular, as shown in 
Section 2, the
operators $Q^\pm_{\lambda,\Lambda}$ can be ``approximated" in $L^\infty$ (the coefficients)  
 by a sequence of operators
with $C^\infty$ coefficients, for which the critical exponent in the radially symmetric case is $p_N^*$. 

The second property is related to the non-monotonicity of the critical 
exponents. 
Notice the following operator's inequalities,
$$
\lambda \Delta\le \mmp\quad\mbox{and}\quad Q^+_{\lambda,\Lambda}\le \mmp,
$$
while for the corresponding critical exponents we have
$$
p^N_*< p^*_+ \quad\mbox{and}\quad \frac{\widetilde N_++2}{\widetilde N_+-2}> p^*_+.$$

We finally observe that all operators of the form ${\cal M}_{s,S}^\pm$ and
$Q^\pm_{s,S}$, with $s,S\in [\lambda,\Lambda]$, have critical exponents in the
interval
$$
\left[\quad\frac{\widetilde N_-+2}{\widetilde N_--2}\,,\quad \frac{\widetilde N_++2}{\widetilde N_+-2}\quad\right].
$$
We conjecture that in the class of operators defined by \equ{ineg_extremale},
the critical exponents are all in the same interval, that is, 
the operators $Q^{\pm}_{\lambda,\Lambda}$ are 
extremal for critical exponents.

\medskip

This is  article is organized in two sections. In section 2 we discuss the case of the operator $Q^+_{\lambda,\Lambda}$. We 
first consider the radial case analyzing, in analogy with the case of the
Laplacian, 
 the critical exponent 
and proving that in the subcritical case the positive solution is 
non-degenerate. Then, by linearization, we show that this solution is non-degenerate in the
space of functions not necessarily symmetric and we apply a perturbation argument via degree theory. In Section 3 we consider the case of the operator
$\mmp$. The situation here is somehow simpler since we know that in the ball all positive solutions are radially symmetric by using a moving plane argument, 
and then, in the subcritical case,  the radial solution is isolated.  We conclude by a 
homotopy invariance with respect to the ellipticity constant, as 
in \cite{busca-esteban-quaas} and then by using a 
perturbation argument.

\section{The extremal operator $Q^+_{\lambda,\Lambda}$.}
\setcounter{equation}{0}

In this section we analyze the equation \equ{1}-\equ{2} in the case
of the extremal operator $Q^+_{\lambda,\Lambda}$, with $0<\lambda<\Lambda$. 
This is a 
uniformly elliptic operator, whose coefficients 
have a discontinuity at the origin. This  feature is what makes this operator 
interesting. We need to make precise the very notion of a solution for
equation \equ{1}-\equ{2}.

We observe that the operator ${Q}^0$ corresponds to the second derivative with respect 
to the radial coordinate $r=|x|$, that is
$$
{ Q}^0 u=\frac{\partial^2 u}{\partial r^2}.
$$
 Because of this observation, we see that 
the analysis of the radial case is very
simple, it amounts to change the notion of the dimension taking into account
$\lambda$ and $\Lambda$, and to perform a phase plane analysis. In fact, we easily see that if
$u$  is a solution of 
\be \label{eradial}
\Lambda u''+\lambda \frac{N-1}{r}u'+u^p=0,\quad u'(0)=u(1)=0.
\ee
then $u(x)=u(|x|)$ is a solution of \equ{1}-\equ{2}.
Defining 
$$
\widetilde N_+=\frac{\lambda}{\Lambda}(N-1)+1 \quad\mbox{and}
\quad v(r)=(\Lambda)^{1/(1-p)}u(r)
$$ 
we
see that $v$ satisfies
\be\label{radial}
v''+\frac{\widetilde N_+-1}{r}v'+v^p=0,\quad v'(0)=v(1)=0.
\ee
The following theorem gives the critical exponent for equation \equ{radial}
and the non-degeneracy property of its solutions in the subcritical case.
\begin{theorem}\label{radialteo}
Equation \equ{radial} does not have a positive solution if
$$
p\ge (\widetilde N_++2)/(\widetilde N_+-2),
$$
and it
possesses exactly one positive solution if
$$
1<p<(\widetilde N_++2)/(\widetilde N_+-2)\quad\mbox{and}\quad \widetilde N_+>2,
$$
or
$$
1<p\quad\mbox{and}\quad 1< \widetilde N_+\le 2.
$$
Moreover, if $v$ is a solution of \equ{radial}, then the linearized equation
\be\label{linear}
h''+\frac{\widetilde N_+-1}{r}h'+pv^{p-1}h=0,\quad h'(0)=h(1)=0,
\ee 
has no non-trivial solution, that is, $0$ is not in the spectrum of 
the linearized operator.
\end{theorem}
{\bf Proof.} The criticality of the number $(\widetilde N_++2)/(\widetilde N_+-2)$
can be proved in a way similar to the case of the Laplacian, by using the
Emden-Fowler transformation. Let $v$ be a positive solution of equation \equ{radial},
then we have that
$$
v_\gamma(r)=\gamma v(\gamma^{(p-1)/2}r),
$$
also satisfies the equation in \equ{radial} together with boundary conditions $v_\gamma'(0)=0$ and
$v_\gamma(\gamma^{(1-p)/2})=0$, for all positive $\gamma$. 
From here we see that the function
$$
h_1(r)=\frac{\partial v_\gamma}{\partial\gamma}(r)|_{\gamma=1},
$$
satisfies \equ{linear} and $h_1'(0)=0$, $h_1(0)>0$.

Assume that $h_2$ is a second solution, linearly independent of $h_1$. Then
necessarily we have that $h_2'(r)$ stays away from zero, for $r$ near $0$, since the 
contrary implies that $h_1$ and $h_2$ are linearly  dependent.  Now, given any
solution $h$ of \equ{linear}, we have
$h=c_1 h_1+c_2 h_2$. But then $c_2=0$ since $h'(0)=0$ and
$c_1=0$ since $h_1(1)<0$, proving that $h\equiv 0$.
%
 $\Box $

\medskip

Continuing with our analysis, we observe that
%
since the operator ${ Q}^0$, and then also ${ Q^+_{\lambda,\Lambda}}$, has discontinuous 
coefficients, we should start making precise the notion of solution for the equation
\begin{eqnarray}
{ Q^+_{\lambda,\Lambda}} u&=&f\quad\mbox{in}\quad \Omega,\label{EQ1}\\
u&=0&\quad\mbox{on}\quad \partial\Omega. \label{EQ2}
\end{eqnarray}
For notational simplicity, in the rest of this section we simply write $Q$
for ${ Q^+_{\lambda,\Lambda}}$, since no confusion will arise.

Given $i,j$ we consider a sequence of $C^\infty$ functions $a_{i,j}^n$ so that
$a_{i,j}^n(x)=x_ix_j/|x|^2$ for all $|x|ŋ\ge 1/n$ and $|a_{i,j}^n(x)|\le 1$
for $|x|ŋ\le 1/n$. For example, we may consider 
a cut-off function $\eta$ so that
$\eta(r)=0$ if $r<1/2$ and $\eta(r)=1$ if $r\ge 1$, and then define
\begin{equation}
a_{i,j}^n(x)=\eta(nr)\frac{x_ix_j}{|x|^2}\label{aprox}.
\end{equation}
Then we define the operators
$$
{ Q}^nu=\lambda \Delta u+(\Lambda-\lambda)\sum_{i=1}^n\sum_{j=1}^na_{i,j}^n
\frac{\partial^2 u}{\partial x_i \partial x_j}.
$$

We assume that the function $f$ is continuous in $\overline \Omega$,
then the problem
\begin{eqnarray}
{ Q}^n u&=&f\quad\mbox{in}\quad \Omega,\\
u&=&0\quad\mbox{on}\quad \partial\Omega,
\end{eqnarray}
possesses a unique smooth solution $u_n$. Moreover, since the coefficients are 
$C^\infty$ functions it is well known the existence of a Green function
$G_n:\Omega\times\Omega \to \R$ allowing to represent this solution as
$$
u_n(x)=\int_\Omega G_n(x,y)f(y)dy.
$$
It follows from Alexandrof-Bakelman-Pucci's estimate that the sequence
$\{G_n(x,\cdot)\}$ is bounded in $L^{N/(N-1)}(\Omega)$ and hence, up to
subsequence, it has a weak limit in this space. Moreover, since our operator ${ Q}$
is discontinuous just at one point, the origin, the weak limit is unique as
shown by the arguments in Cerutti, Escauriaza and Fabes \cite{cerutti}.
Thus, our problem has a unique Green function $G:\Omega\times\Omega \to \R$, 
such that $G(x,\cdot)$ is in $L^{N/(N-1)}(\Omega),$ and we  define
$$
u(x)=\int_\Omega G(x,y)f(y)dy,
$$ 
as the solution to equation \equ{EQ1} and \equ{EQ2}.

On the other hand, the sequence of 
 solutions $\{u_n\}$
is bounded in $C^\alpha_0(\bar \Omega)$, for $\alpha>0$ as follows from
basic estimates, see Gilbarg and Trudinger \cite{GT} or Cabr\'e Caffarelli 
\cite{cc}. Then, the solution $u$ is actually of class $C^\alpha(\bar \Omega)$ 
and we have
\be
\|u\|_{C^\alpha(\overline \Omega)}\le C\|u\|_{L^\infty(\partial\Omega)}+
C\|f\|_{L^N(\Omega)}.
\label{est}\ee
Here the constant $C$ depends only on the
ellipticity constants, the $L^\infty$ bounds on the coefficients and on the
domain $\Omega$, which we assume to have a regular boundary.

%
%

Let $u_0$ be the unique solution of \equ{eradial}. 
In what follows we 
show that, when $u_0$ is  considered as a function of the $N$ variables, it satisfies
the equation in the sense given above.
\begin{lemma}
The function $u_0(x)=u_0(|x|)$ satisfies the equation 
\begin{eqnarray}
{ Q} u&=&-u^p\quad\mbox{in}\quad B,\label{ball1}\\
u&=&0\quad\mbox{on}\quad \partial B,\label{ball2}
\end{eqnarray}
in the sense just defined above.
\end{lemma}
\noindent
{\bf Proof.} By direct computation we see that, pointwise, we have
$$
{ Q}^n u=-u^p -c_n(x)\quad\mbox{in}\quad B,
$$
where $c_n$ is a function with support in the ball $B(0,1/n)$ and which is
bounded, with
a bound independent of $n$. Then we certainly have
$$
u(x)=\int_B G_n(x,y)(-u^p(y)+c_n(y)) dy.
$$
Taking limits here we conclude. $\Box$

\medskip

\noindent
\bremark The notion of solution defined above is known 
as {\it good solution} and it was introduced by Cerutti, 
Escauriaza and Fabes \cite{cerutti}. In a recent paper by Jensen, Kocan and 
Swiech \cite{jensen}, 
this notion of solutions is shown to be equivalent to $L^p$ viscosity 
solution.
\eremark
\medskip
                                                                                
\noindent
\bremark
 We do not know whether equation \equ{ball1}-\equ{ball2} 
possesses a non-radial solution or not.
\eremark
\medskip

Our existence result is for domains which are close to the unit ball. 
More precisely we
 assume that we have a sequence of domains $\{\Omega_n\}$ 
such that for all $0<r<1<R$ there exists
$n_0\in \N$ such that
$$
B(0,r)\subset \Omega_n\subset B(0,R),\quad \mbox{for all}\,\, n\ge n_0.
$$
We consider a reference bounded domain $D$ such that $\Omega_n\subset D$, for all
$n\in\N$. We may take, for example, $D$ as the ball of radius $2$.

Next we prove a continuity property for the Green functions associated to
the domains $\Omega_n$
\begin{lemma}\label{conv}
Under the conditions given above we have, for every $f\in C(D)$ and
for every $x\in B$,
$$
\lim_{n\to\infty}\int_{\Omega_n}G_{\Omega_n}(x,y)f(y)dy=\int_{B}G_B(x,y)f(y)dy,
$$
where $G_B$ and $G_{\Omega_n}$ are  the Green functions of the unit ball $B$
and of $\Omega_n$, respectively.
\end{lemma}
\noindent
{\bf Proof.}
Let $f\in C(D)$ and let $u_n$ be the solution of the equation
$$
{ Q} u= f,\quad\mbox{in}\quad \Omega_n,\qquad  u=0  \quad\mbox{on}
\quad \partial \Omega_n.
$$
Consider also the function $u$, the solution of the equation
$$
{ Q} u= f,\quad\mbox{in}\quad B_R,\qquad  u=0  \quad\mbox{on}
\quad \partial B_R.
$$
Here we assume $r<1<R$ are close to $1$ and $n$ is large enough so that $B_r
\subset \Omega_n\subset B_R$. Using \equ{est} we find
$$
\|u\|_{C^\alpha(\overline B_R)}\le 
C\|f\|_{L^N(B_R)},
$$
which implies
$$
|u(x)|\le C(R-r)^\alpha \|f\|_{L^N(B_R)},\quad x\in B_R\setminus B_r.
$$
Since $u_n-u$ satisfies ${ Q}(u_n- u)=0$ in $\Omega_n$, by standard estimates,
see \cite{cc}, we obtain
$$
\|u_n-u\|_{L^\infty(\Omega_n)}\le C \|u_n-u\|_{L^\infty(\partial\Omega_n)}\le
C(R-r)^\alpha \|f\|_{L^N(B_R)},
$$
where the constant $C$ here is uniform in $n$.
Similarly if $v$ is the solution of
$$
{ Q} u= f,\quad\mbox{in}\quad B,\qquad  u=0  \quad
\mbox{on}
\quad \partial B,
$$
then we have
$$
\|v-u\|_{L^\infty(B)}\le 
C(R-r)^\alpha \|f\|_{L^N(B)}.
$$
Since for all $x\in B_r$ we have
$$
u_n(x)-v(x)=
\int_{\Omega_n}G_{\Omega_n}(x,y)f(y)- \int_BG_{B}(x,y)f(y),
$$
the result follows.
$\Box$

\medskip

Now we state and prove our main theorem of this section.
\begin{theorem}\label{teoprincipal}
Assume $\tilde N_+>2$ and that $1<p<(\tilde N_++2)/(\tilde N_+-2)$. Then there is $n_0\in\N$ so that for all $n\ge n_0$,
the equation
\begin{eqnarray}\label{teo1}
{ Q} u+u^p&=&0\quad\mbox{in}\quad \Omega_n,\label{EQ11}\\
u&=&0\quad\mbox{on}\quad \partial\Omega_n, \label{EQ21}
\end{eqnarray}
possesses at least one nontrivial solution.
\end{theorem}
\noindent
In order to prove our theorem we will follow some ideas from \cite{dancer}.
We start setting up the functional analytic framework.
We consider the inclusions $i:C_0(\bar D)\to C(\bar B)$ and
$j:C_0^\alpha(\bar D )\to C_0(\bar D)$. Here $C$ stands for continuous functions, $C_0$ for continuous functions vanishing on the boundary and $C_0^\alpha$
for H\"older continuous functions vanishing on the boundary. 
Given  $f\in C(\overline B)$, we let $\bar {\cal L}(f)$ be the unique solution to 
\equ{EQ1}-\equ{EQ2} with $\Omega=B$. Then we extend this solution
to define ${\cal L}(f)$ as
\begin{equation}\label{linearLL}
{\cal L}(f)(x)=\left\{ \begin{array}{ll}
        \bar {\cal L}(f)(x)   &\mathrm{if}\quad x\in \bar B,
        \\ 0  &\mathrm{if}\quad x\in \bar  D\setminus B.
        \end{array}\right.
\end{equation}
Thus, the operator ${\cal L}:C(\overline B)\to
C^\alpha_0(\overline D)$ is well defined as a linear bounded operator.
Next we define the nonlinear operator ${\cal F}:C_0(\overline D)\to C_0(\overline D)$ as
\begin{equation}
{\cal F}(u)=-j\circ{\cal L}(i(u^p)).\label{calf}
\end{equation}
If we consider $N: C_0(\overline D)\to C_0(\overline D)$, the Nemitsky operator
defined as $N(u)=u^p$, we easily see that $N$ is of class $C^1$  and 
$N'(u)(h)=pu^{p-1}h$ for $h\in C_0(\overline D)$.
Thus, the operator ${\cal F}$ is compact and of class $C^1$.

In the definition of the operators $i$ and ${\cal L}$ we can replace $B$ by
$\Omega_n$ and we obtain $i_n$ and ${\cal L}_n$. The operator
${\cal F}_n$ is then defined as ${\cal F}_n(u)=-j\circ{\cal L}_n(i_n(u^p))$,
and naturally it is compact and of class $C^1$ as an operator in 
$C_0(\overline D)$.

\medskip

Before giving the proof of Theorem \ref{teoprincipal} we need still another 
preliminary result. It guarantees that $u_0$ is isolated not only in the space
of radial functions, as shown in Theorem \ref{radialteo}, but also in the space of all functions of class $C^\alpha$ in the ball $B=B(0,1)$.
\begin{proposition}\label{nondeg}
Under the hypotheses of Theorem \ref{teoprincipal}, the linear equation
\begin{eqnarray}\label{linearL}
{ Q} h+pu_0^{p-1}h&=&0\quad\mbox{in}\quad B,\label{EQli1}\\
h&=&0\quad\mbox{on}\quad \partial B, \label{EQli2}
\end{eqnarray}
has only the trivial solution $h\equiv 0$ in $C_0^\alpha(B)$.
\end{proposition}
\noindent
{\bf Proof.}
Our proof uses a standard argument by Smoller and Wasserman \cite{smollerwasserman}. 
We let $\{\phi_k(\theta)\}_{k=0}^\infty$, $\theta\in S^{N-1}$, be the 
eigenfunctions of the Laplacian in $S^{N-1}$, whose eigenvalues are
$$
\lambda_k=-k(k+N-2).
$$
Assume $h$ is a solution of \equ{EQli1}-\equ{EQli2}. Then we 
consider the approximation operator ${ Q}^n$ as given by \equ{aprox}
and we solve the equation
\be\label{pp}
{ Q}^nh_n=-pu_0^{p-1}h \quad\mbox{in}\quad B,\quad  h_n=0,\quad\mbox{on}\quad
\partial B.
\ee
We observe that if  $u\in C^2$ and  $u=u(r,\theta)$ with 
$ \theta\in S^{N-1}$ then
for the operator   ${ Q}^n$ we can write
$$
{ Q}^nu= 
\{\lambda+(\Lambda-\lambda)\eta(nr)\}u''
+\lambda\frac{N-1}{r}u'+\frac{\lambda}{r^2}\Delta_\theta u,
$$
where $'$ denotes derivative with respect to $r$ and $\Delta_\theta$ is the
Laplacian on the sphere.
Let us define 
$$
a_k^n(r)=\int_{S^{N-1}} h_n\phi_k\, d\theta\quad\mbox{
and }
a_k(r)=\int_{S^{N-1}} h\phi_k\, d\theta.
$$
By the $C^{\alpha}$ convergence of $h_n$ to $h$ we see that $a_k^n$ converges 
uniformly in $[0,1]$ to
$a_k$, for all $k\ge 0$.
Next we multiply equation  \equ{pp} by $\phi_k$ and integrate over $S^{N-1}$
to obtain
$$
\{\lambda+(\Lambda-\lambda)\eta(nr)\}(a_k^n)''
+\lambda\frac{N-1}{r}(a_k^n)'+\lambda\lambda_k\frac{a_k^n}{r^2}= -pu_0^{p-1}a_k.
$$
From here we can prove that for $r\in (0,1]$ 
the convergence of $a_k^n$ to $a_k$ is even
$C^2$ and $a_k$ satisfies
\begin{equation}\label{dd}
a_k''+\frac{\tilde N-1}{r}a_k'+\frac{\lambda\lambda_k}{\Lambda}\frac{a_k}{r^2}
+\frac{p}{\Lambda}u_0^{p-1}a_k=0.
\end{equation}
Moreover, we can prove that $a_k(0)=0$ and that $r^{\tilde N_+-1}a_k'(r)$ 
is bounded as $r\to 0$. To prove the last statement we use that $\tilde N_+>2$.
On the other hand we have that $w=u_0'$ satisfies the equation
\begin{equation}\label{dd1}
w''+\frac{\tilde N-1}{r}w'-\frac{\tilde N-1}{r^2}w+
\frac{p}{\Lambda}u_0^{p-1}w=0.
\end{equation}
Multiplying equation \equ{dd} by $w$ and equation \equ{dd1} by $a_k$, integrating between
$0$ and the first zero of $a_k$, and subtracting we can prove that
$a_k\equiv 0$ for all $k\ge 1$, see \cite{smollerwasserman}. For $k=0$
we use Theorem \ref{radialteo} to prove that also $a_0\equiv 0$.$\Box$

\medskip

We finally complete the proof of our main theorem in this section.

\noindent
{\bf Proof of Theorem \ref{teoprincipal}.}
Since $I-{\cal F}'(u_0)$ has trivial kernel, as we proved in 
Proposition \ref{nondeg}, 
there is $\delta>0$ such that $u-{\cal F}(u)\not =0$
for all $u\in\partial {\cal B}$, where
${\cal B}=\{u\in C_0(\bar D)\,/\, \|u-u_0\|_{C_0(\bar D)}<\delta\}$.
Moreover, the Leray-Schauder degree of $I-{\cal F}$ is well defined in ${\cal B}$ and 
 ${\rm deg}(I-{\cal F}, {\cal B},0)=1$ or $-1$.
To
finish the proof we just need to prove that there exists
$n_0\in\N$ so that 
\begin{equation} \label{B}
u\not = t{\cal F}(u)+(1-t){\cal F}_n(u), \quad \mbox{for all}\quad t\in[0,1],\,
u\in\partial {\cal B},
\end{equation}
since
 this implies that
$$
{\rm deg}(I-{\cal F}_n, {\cal B},0)=
{\rm deg}(I-{\cal F}, {\cal B},0)\not =0.
$$
Let us assume that \equ{B} is not true. Then there exist sequences $\{t_n\}\subset [0,1]$
and $\{u_n\}\subset \partial {\cal B}$ such that
$$
u_n=t_n{\cal F}(u_n)+(1-t_n){\cal F}_n(u_n).
$$
We may assume that 
$u_n\to\bar u$ uniformly in $\bar D$, up to a subsequence, as a 
consequence of the compactness of the inclusion $j$. 
We may also assume $t_n\to \bar t$.
Then we see that we will get a contradiction if we prove that
\begin{equation} \label{BB}
\lim_{n\to \infty} {\cal F}_n(u_n)={\cal F}(\bar u).
\end{equation}
From the definition of ${\cal F}_n$ we see that if $v_n={\cal F}_n(u_n)$ then
$$
v_n(x)=-\int_{\Omega_n}G_{\Omega_n}(x,y)u_n^p(y)dy,
$$
where $G_{\Omega_n}$ is the Green function of $\Omega_n$.
Since $u_n$ is uniformly convergent to $\bar u$ in $\bar D$,  
to complete the proof we just use Lemma \ref{conv}.
$\Box$

\medskip
\noindent
\bremark Our Theorem \ref{teoprincipal} is concerned with the operator
$Q^+_{\lambda,\Lambda}$, for $0<\lambda\le\Lambda$. A completely analogous theorem
can be proved for the operator $Q^-_{\lambda,\Lambda}$. Naturally, our hypotheses have to be changed to: $\tilde N_->2$ and $1<p<(\tilde N_-+2)/(\tilde N_--2)$.
\eremark

\setcounter{equation}{0}
\section{The extremal operator $\mmp$.}
                                                                                
In this section we analyze equation \equ{1}-\equ{2} in the case $L= M^+_{\lambda,\Lambda}$, the extremal Pucci's operator for $0<\lambda\le \Lambda$. This is an
elliptic operator which is nonlinear, but homogenous of degree 1.
In comparison with ${ Q^+_{\lambda,\Lambda}}$, for the operator
$\mmp$ there is a good regularity
theory that guarantees that the solutions to equation \equ{1}-\equ{2} are classical solutions (see \cite{cc}).

We start this section recalling  the existence of positive solutions
 for \equ{1}-\equ{2}, when  $\Omega$ 
is a ball as was proved by  Felmer and Quaas 
in \cite{felmerquaas1} and \cite{felmerquaas}.
\begin{theorem} \label{teo11} Suppose  $\tilde N_+> 2$. Then there exist a number 
$p^*_+$ such that 
$$\frac{N+2}{N-2}<p^*_+< \frac{\widetilde N_++2}{\widetilde N_+-2}\,.$$ 
with the property that
 if $1<p<p^+_*$  then \equ{1}-\equ{2} has a nontrivial radially  symmetric $C^2$ solution.
Moreover, if $p \geq p^+_*$ then \equ{1}-\equ{2} does not have a solution.
\end{theorem}
The number $p^+_*$ is called critical exponents for the operator
$M^+_{\lambda,\Lambda}$. 
  We notice that  $\lambda=\Lambda$ implies $p^+_*=p_*^N$.

Let $\Omega_n$ be a sequence of domains satisfying the conditions given
in the previous section. Now we present the main theorem of this section
\begin{theorem} \label{teo2}
Assume $\tilde N_+ >2$ and that $1< p<p_+^*$. Then there is $n_0\in\N$ so that for all $n\ge n_0$,
the equation
\begin{eqnarray}
\mmp (D^2 u)+u^{ p}&=&0\quad\mbox{in}\quad \Omega_n,\label{EQp11}\\
u&=&0\quad\mbox{on}\quad \partial\Omega_n, \label{EQp21}
\end{eqnarray}
possesses at least one positive solution.
\end{theorem}
\bremark 
It can be seen  that equation
\equ{1}-\equ{2}, when $\Omega$ is a ball, has only one 
positive radially symmetric solution. In fact,
by using a classical moving planes technique (see \cite{B}), all positive solutions of \equ{1}-\equ{2} are radially symmetric 
(for detail see \cite{sirakov}). 
This uniqueness property  is crucial in our analysis, since it allows us
to use a degree theory approach as in the previous section, avoiding the study
of the linearized equation in order to obtain non-degeneracy of the radially symmetric solution.
\eremark
\bremark 
a) Our Theorem \ref{teo2} is concerned with the operator
$M^+_{\lambda,\Lambda}$, for $0<\lambda\le\Lambda$. A completely analogous theorem
can be proved for the operator $M^-_{\lambda,\Lambda}$. With the 
natural change in the hypothesis  to $1<p<p_-^*$.

\noindent b) For other existence result concerning Pucci operator we refer the reader to \cite{felmerquaas2} and \cite{quaas1}.
\eremark
\noindent {\bf Proof.} Let $q:[\lambda, \Lambda] \to \RR$ be a continuous function such that 
$q(s) < p^+_*(s)$ for all $s \in [\lambda, \Lambda]$ and $q(\lambda)= p$. Here $p^+_*(s)$ is the critical exponent for the operator ${\cal M}_{s,\Lambda}^+$, which is a continuous function of $s$.
Next we consider $s\in [\lambda, \Lambda]$ and we define 
$\bar{\cal L}^s(f)$ as 
the unique solution to   
 \begin{eqnarray}
 {\cal M}_{s,\Lambda}^+ (D^2 u)&=&f\quad\mbox{in}\quad B,\label{EQP1}\\
          u&=0&\quad\mbox{on}\quad \partial B, \label{EQP2} 
                              \end{eqnarray}
for $f \in C(\bar B)$. This operator is well defined in $C(\bar B)$, 
with
values in $C^\alpha_0(\bar B)$, and it is positive, that is 
if $f(x)\le 0$ in $B$ 
then $u\ge 0$ in $B$. These follow from existence and regularity theory
for fully nonlinear operators and the maximum principle for 
${\cal M}_{s,\Lambda}^+$. See 
the monography by  Caffarelli and Cabr\'e  \cite{cc} and the work by
Bardi and Da Lio \cite{bardidalio}.
Now we extend $\bar { \cal L}^s$ as in \equ{linearLL} and  define the operator 
${\cal F}(s,\cdot)$ using  \equ{calf}, 
 with $p=q(s)$.

Let $u^s_0$ be the unique solution of \equ{1}-\equ{2} in $B$ for $p=q(s)$ given by Theorem \ref{teo11}. 
We consider
$$M=\sup_{s \in [\lambda, \Lambda]}\| u^s_0\|_{L^\infty(B)}
\quad\mbox{and}\quad m=\inf_{s \in [\lambda, \Lambda]}\| u^s_0\|_{L^\infty(B)}
>0,
$$
and let $\delta>0$ be such that $m-\delta>0$. We define the set
 $${\cal B}=\{u\in C_0(\bar D)\,/\, u(x)\ge 0 \mbox{ in } D,\,\, 
m-\delta <\|u\|_{C_0(\bar D)}<M+1\},
$$
and we observe, as follows 
by the uniqueness of solutions in the ball and the strong maximum principle for
${\cal M}_{s,\Lambda}^+$, that
$u-{ F}(s, u)\not =0$ for all $u\in\partial {\cal B}$ and 
$s \in [\lambda, \Lambda]$. 

We notice that  ${\rm deg}(I-{\cal F}(\Lambda,\cdot), {\cal B},0)\not =0$, since
in this case we are dealing with the Laplacian. Then by 
invariance under homotopy of the degree we obtain
$$
{\rm deg}(I-{\cal F}{(\lambda,\cdot)}, {\cal B},0)=
{\rm deg}(I-{\cal F}{(\Lambda,\cdot)}, {\cal B},0)\not =0.
$$
Next we define ${\cal F}_n$ perturbing the domain as in the previous section, 
keeping $s=\lambda$ and $p=q(\lambda)$. 
To
finish the proof we just need to prove that there exists
$n_0\in\N$ so that for all $n \geq n_0$ 
\begin{equation} \label{Bp}
u\not = t{\cal F}( \lambda,u)+(1-t){\cal F}_n(\lambda,u), \quad \mbox{for all}\quad t\in[0,1],\,
u\in\partial {\cal B},
\end{equation}
since
 this implies that
$$
{\rm deg}(I-{\cal F}_n(\lambda,\cdot), {\cal B},0)=
{\rm deg}(I-{\cal F}(\lambda,\cdot), {\cal B},0)\not =0.
$$
Let us assume that \equ{Bp} is not true. Then there exist sequences $\{t_n\}\subset [0,1]$
and $\{u_n\}\subset \partial {\cal B}$ such that
$$
u_n=t_n{\cal F}( \lambda,u_n)+(1-t_n){\cal F}_n(\lambda,u_n).
$$
We may assume that 
$t_n\to \bar t$ and that
$u_n\to\bar u$ uniformly in $\bar D$, up to a subsequence, as a 
consequence of the compactness of the inclusion $j$. 
We will get a contradiction if we prove that
\begin{equation} \label{BBp}
\lim_{n\to \infty} {\cal F}_n( \lambda,u_n)(x)={\cal F}( \lambda,\bar u)(x)\quad \mbox{for all } x \in D.
\end{equation}
We first notice that ${\cal F}_n( \lambda,u_n):=v_n \to v$ in $C_0(D)$. 
If $x \in D \setminus \bar B$, then for $n$ large  ${\cal F}_n( \lambda,u_n)(x)=0$, and  so $v(x)=0$, 
therefore by continuity of $v$ we also have $v=0$ in $D \setminus  B$.
On the other hand 
$$\mmp(D^2 v_n)=-u^{ p}_n\quad\mbox{in}\quad B(0, r)$$ for $n$ large.
Passing to the limit in the viscosity sense we get 
$\mmp(D^2 v)=-\bar u^{ p}$ in $B$. 
Since 
$v=0$ on $ \partial B $ then, 
by the definition of ${\cal F}$, we obtain ${\cal F}( \lambda, \bar u)=v$, 
concluding the proof. $\Box$
\medskip

\noindent

{\bf Acknowledgements}
The second author was  partially supported by Fondecyt Grant 
\# 1030929 and FONDAP de Matem\'aticas Aplicadas. The third author 
was partially supported by Fondecyt Grant \# 1040794. This work was 
partially supported by ECOS Grant \# C02E08.

E-mail addresses: esteban@ceremade.dauphine.fr, pfelmer@dim.uchile.cl,\\
 alexander.quaas@usm.cl

\end{document}